%% file: article.tex
\documentclass[titlepage,14pt]{extarticle}
\usepackage{amsfonts,amsmath,amssymb,indentfirst}

\usepackage{geometry}
\newtheorem{theorem}{Theorem}
\newtheorem{lemma}{Lemma}
\newtheorem{corollary}{Corollary}
\def\proof{{\bf Proof. }}

\def\inbr#1{ \{ #1 \} }
\def\cube#1{[0,1]^{#1}}
\newcommand{\Le}{\mathop{\rm Le}\nolimits}
\newcommand{\sign}{\mathop{\rm sign}\nolimits}
\def\C{{\mathbb C}}
\def\Z{{\mathbb Z}}
\def\N{{\mathbb N}}
\def\Q{{\mathbb Q}}
\def\F{{\mathcal{F}}}
\newcommand{\nomer}{{No. }}

\allowdisplaybreaks[4]

\title{Properties of Coefficients of Certain Linear Forms
in Generalized Polylogarithms.}
\author{Zlobin~S.A.}

\begin{document}

\maketitle

{\flushleft UDK 511.36}

{\it Abstract. }
We study properties of coefficients of a linear form,
originating from a multiple integral.
As a corollary, we prove Vasilyev's conjecture,
connected with the problem of irrationality of the Riemann zeta function
at odd integers.


\section{Introduction}
\input intro.tex

\section{Elementary Sums}
\input ElemSums.tex

\section
{Denominators of linear form coefficients}
\label{ArithmSection}
\input Arithm.tex

\section
{Estimate of linear form coefficients}
\label{UpperBoundSection}
\input UpperBound.tex


\newcommand{\namefont}{\scshape}
\newcommand{\titlefont}{\itshape}

\end {document}

%% file: intro.tex
Define a generalized polylogarithm by the series
$$
\Le_{\vec{s}}(z) = \sum_{n_1 \ge n_2 \ge \dots \ge n_l \ge 1 }
\frac{z^{n_1}}{n_1^{s_1} n_2^{s_2} \dots n_l^{s_l}}
$$
for a vector $\vec{s}=(s_1,\dots,s_l)$, $s_j \in \N$. 
This series converges when $|z|<1$.
In the paper \cite{zl5},
in connection with approximations of values of
generalized polylogarithms, it has been proved the following result.
Let $a_i$, $b_i$, $c_j$ be integers, satisfying the certain inequalities.
Then, the 
identity
\begin{equation}
\label{linform}
S(z) = \int_{\cube{m}} \frac{ \prod_{i=1}^m
x_{i}^{a_{i}-1} (1-x_{i})^{b_{i}-a_{i}-1}}
{\prod_{j=1}^l (1-z x_1 x_2 \dots x_{r_j})^{c_{j}} } dx_1 dx_2
\dots dx_{m} =
\sum_{\vec{s}} P_{\vec{s}}(z^{-1}) \Le_{\vec{s}}(z),
\end{equation}
holds, where $0=r_0 < r_1 < r_2 < \dots < r_l = m$ and
$P_{\vec{s}}$ are polynomials with rational coefficients.
This representation is unique 
because of the linear independence of
$\Le_{\vec{s}}(z)$ with different indices
over $\C(z)$ (see \cite[Corollary 1]{zl5}).

It is important in arithmetical applications
to have estimate for absolute values and a common denominator of
coefficients of the polynomials $P_{\vec{s}}$
depending on the parameters $a_i$, $b_i$, $c_j$,
as well as on the dimension $m$ of the integral.
This is the main aim of this paper.

One of possible applications of integrals of type $S(z)$
is related to the problem of the irrationality of the Riemann
zeta function  $\zeta(k)$ at odd integers $k=3,5,7,\dots$.
In \cite{vasiliev}, D.V.~Vasilyev considered the integrals
$$
V_{m,n} = \int_{\cube{m}} \frac{\prod_{i=1}^m x_i^{n}
(1-x_i)^{n} } {\left( 1-x_1(1-x_2(\cdots-x_{m-1}(1-x_m)\cdots
\right)^{n+1}} dx_1 dx_2 \dots dx_m.
$$
He conjectered that
\begin{equation}
\label{VasConj}
V_{2l+1,n}=A_0 + \sum_{j=1}^{l}
A_j \zeta(2j+1), \quad D_n^{2l+1} A_{j} \in \Z,
\end{equation}
where $D_n$ is the least common multiple of 1, 2, \dots, $n$.
The integral $V_{3,n}$ is equal to 
the integral, which was used by F.~Beukers for the
proof of the irrationality of $\zeta(3)$ (see \cite{beukers}).
The equality (\ref{VasConj}) holds for it.
Vasilyev proved (\ref{VasConj}) for $m=5$.
Later W.V.~Zudilin (\cite{zudilin}) showed (\ref{VasConj})
with the weaker inclusion
$D_n^{2l+2} \Phi_n^{-1} A_j \in \Z$,
where $\Phi_n$ is the product of prime numbers $p<n$
for which $2/3 \le \{ n/p \} < 1$
($\{ \cdot \}$ denotes the fractional part of a number).
The validity of $D_n^{2l+1} A_j \in \Z$ was proved
by C.~Krattenthaler and T.~Rivoal
(\cite[Th\'eor\`eme 1]{rivkrat}).
Their proof is technically complicated.
In this paper we prove (\ref{VasConj})
using the following representation $V_{m,n}$ in the form (\ref{linform})
(see \cite[Corollary 2]{zl1}):
\begin{equation}
\label{VEqualSOdd}
V_{2l+1,n}= \int_{\cube{2l+1}} \frac{ \prod_{i=1}^{2l+1}
x_{i}^n (1-x_{i})^n dx_1 dx_2 \dots
dx_{2l+1}} { \prod_{j=1}^l (1-x_1 \dots
x_{2j} )^{n+1} (1-x_1 x_2 \dots x_{2l} x_{2l+1})^{n+1} }.
\end{equation}

We prove theorems \ref{ArithmTh1} and \ref{UpperBoundTh}
in sections \ref{ArithmSection} and \ref{UpperBoundSection} of this article. 
They give
the estimate on the common denominator and the values of the 
coefficients $P_{\vec{s}}$ in (\ref{linform})
under the certain conditions.

%% file: ElemSums.tex
We call a sum of type
\begin{equation}
\label{simplesum}
\sum_{n_1 \ge n_2 \ge \dots \ge n_l \ge 1}
z^{n_1-1} \prod_{j=1}^l \frac{1}{(n_j+p_j)^{u_j}},
\quad p_j \in \{0,1,2,\dots\}, \quad u_j \in \N,
\end{equation}
{\it elementary}.
From \cite[Theorem 1]{zl5} it follows that this
sum can be expressed in the form (\ref{linform}).

In what follows, for any vector $\vec{s}=(s_1, s_2, \dots, s_l)$ we use
the notation $w(\vec{s})=s_1+s_2+\dots+s_l$. 
The {\it height} of the polynomial is the maximum of the absolute values of its
coefficients.

\begin{lemma}
\label{shifted}
Let $P=\max\limits_{1 \le j \le l} p_j$.
Then, for the sum (\ref{simplesum}), the heights of the polynomials
$P_{\vec{s}}$ do not exceed
\begin{equation}
\label{ShiftedUpperBound}
\max( l! \cdot (w(\vec{u}) 2^{w(\vec{u})})^{l-1} P^l, 1);
\end{equation}
moreover,
$D_{P}^{w(\vec{u})-w(\vec{s})} P_{\vec{s}}(z) \in \Z[z]$.
\end{lemma}
\proof
We use the following notation:
$r_0=0$, $r_j = u_1 + u_2 + \dots + u_j$, $m = r_l = w(\vec{u})$. 
By \cite[Lemma 2]{zl5}
it is possible to write expression (\ref{simplesum})
as the integral
$$
I(p_1, p_2, \dots, p_l) = \int_{\cube{m}} \frac{\prod_{j=1}^l
(x_{r_{j-1}+1} x_{r_{j-1}+2}\dots x_{r_j})^{p_j}}
{\prod_{j=1}^l (1-z x_1 x_2 \dots x_{r_j}) }
dx_1 dx_2 \dots dx_{m}.
$$

We prove Lemma \ref{shifted} by induction
on the vector $(l,p_1+p_2+\dots+p_l)$.
We order vectors $(l,k)$ in lexicographic ally, i.e.
$$
(l_1, k_1 ) < (l_2, k_2) \Leftrightarrow
l_1 < l_2 \mbox{ or } l_1=l_2 \mbox{ and } k_1<k_2.
$$
The statement, which is proved by induction, is
a little stricter than the 
statement of the lemma: the heights of $P_{\vec{s}}(z)$ do not exceed
$$
\max \left( \sum_{j=1}^l p_j \cdot (l-1)! \cdot (m 2^m P)^{l-1}, 1 \right).
$$
This estimate is really more precise than (\ref{ShiftedUpperBound})
since $\sum_{j=1}^l p_j \le l \cdot P$.
The induction base ($p_1=p_2=\dots=p_l=0$) follows from (\ref{simplesum}):
$I(0,0,\dots,$ $0) = z^{-1} \Le_{u_1, u_2, \dots, u_l}(z)$.

Let $p_h>0$ for some $h>1$.
From the equality
\begin{eqnarray*}
( x_{r_{h-1}+1} x_{r_{h-1}+2}\dots x_{r_h})^{p_h} & = &
(x_{r_{h-1}+1} x_{r_{h-1}+2}\dots x_{r_h})^{p_h-1} \\
&& + (x_{r_{h-1}+1} x_{r_{h-1}+2}\dots x_{r_h})^{p_h}
( 1 - z x_1 x_2\dots x_{r_{h-1}} )\\
&& - (x_{r_{h-1}+1} x_{r_{h-1}+2}\dots x_{r_h})^{p_h-1}
( 1 - z x_1 x_2\dots x_{r_{h}} )
\end{eqnarray*}
it follows that
\begin{align*}
I(p_1, p_2, \dots,p_h, \dots, p_l) & =
I(p_1, p_2, \dots,p_h-1, \dots, p_l) \\
& \; + \int_{\cube{m}} \frac{\prod_{j=1}^l
(x_{r_{j-1}+1} x_{r_{j-1}+2}\dots x_{r_j})^{p_j}}
{\prod_{{j=1} \atop j \ne h-1}^l (1-z x_1 x_2 \dots x_{r_j}) }
dx_1 dx_2 \dots dx_{m} \\
& \; - \int_{\cube{m}} \frac{\prod_{j=1}^l
(x_{r_{j-1}+1} x_{r_{j-1}+2}\dots x_{r_j})^{p_j'}}
{\prod_{{j=1} \atop j \ne h}^l (1-z x_1 x_2 \dots x_{r_j}) }
dx_1 dx_2 \dots dx_{m},
\end{align*}
where $p_j' = p_j$ for $j \ne h$ and $p_h' = p_h - 1$.
By \cite[Lemma 2]{zl5} we write this equality as
\begin{align}
I(p_1, p_2, \dots,p_h, & \dots, p_l) \nonumber \\
\label{ShiftedCase2Int}
& = I(p_1, p_2, \dots,p_h-1, \dots, p_l) \\
& +
\sum_{n_1 \ge n_2 \ge \dots \ge n_{l-1} \ge 1}
z^{n_1-1}
\prod_{j=1}^{h-2} \frac{1}{(n_j+p_j)^{u_j}} \nonumber \\*
\label{ShiftedCase2Sum1}
& \quad \times
\frac{1}{(n_{h-1}+p_{h-1})^{u_{h-1}} (n_{h-1}+p_{h})^{u_{h}}} \cdot
\prod_{j=h}^{l-1} \frac{1}{(n_j+p_{j+1})^{u_{j+1}}} \\
& -
\sum_{n_1 \ge n_2 \ge \dots \ge n_{l-1} \ge 1}
z^{n_1-1}
\prod_{j=1}^{h-1} \frac{1}{(n_j+p_j)^{u_j}} \nonumber \\
\label{ShiftedCase2Sum2}
& \quad \times
\frac{1}{(n_{h}+p_{h}-1)^{u_{h}} (n_{h}+p_{h+1})^{u_{h+1}}} \cdot
\prod_{j=h+1}^{l-1} \frac{1}{(n_j+p_{j+1})^{u_{j+1}}}
\end{align}
If $h=l$, the subtracted sum reads as
$$
\frac{1}{p_l^{u_l}}
\sum_{n_1 \ge n_2 \ge \dots \ge n_{l-1} \ge 1}
z^{n_1-1}
\prod_{j=1}^{l-1} \frac{1}{(n_j+p_j)^{u_j}}
$$

Now we consider in detail the sum (\ref{ShiftedCase2Sum1}).
If $p_{h-1}=p_h$, then
$$
\frac{1}{(n_{h-1}+p_{h-1})^{u_{h-1}} (n_{h-1}+p_{h})^{u_{h}}}=
\frac{1}{(n_{h-1}+p_{h-1})^{u_{h-1}+u_{h}}},
$$
i.e. the sum (\ref{ShiftedCase2Sum1}) 
is elementary and it is possible to apply the induction hypothesis to it.
In this case the heights of
polynomials $P_{\vec{t}}(z)$ in the sum decomposition
(into a linear form) do not exceed
$$
(l-1)! \cdot (m 2^m)^{l-2} P^{l-1},
$$
and the common denominator of
coefficients of $P_{\vec{t}}(z)$ divides $D_P^{m-w(\vec{t})}$.
If $p_{h-1} \ne p_h$, then we take the following 
partial fraction decomposition:
$$
\frac{1}{(n_{h-1}+p_{h-1})^{u_{h-1}} (n_{h-1}+p_{h})^{u_{h}}}=
\sum_{k=1}^{u_{h-1}} \frac{A_k}{(n_{h-1}+p_{h-1})^k}+
\sum_{k=1}^{u_{h}} \frac{B_k}{(n_{h-1}+p_{h})^k},
$$
$$
A_k=(-1)^{u_{h-1}-k} \binom{u_{h-1}+u_h-k-1}{u_{h-1}-k} \frac{1}
{(p_h-p_{h-1})^{u_{h-1}+u_h-k}},
$$
$$
B_k=(-1)^{u_{h}-k} \binom{u_{h-1}+u_h-k-1}{u_{h}-k} \frac{1}
{(p_{h-1}-p_h)^{u_{h-1}+u_h-k}}.
$$
Substituting this equality into (\ref{ShiftedCase2Sum1}),
we write (\ref{ShiftedCase2Sum1}) as 
the sum of $u_{h-1}+u_h$ elementary sums
(with coefficients $A_k$ and $B_k$). 
We can apply the induction hypothesis for each of them.
Consider one of them,
$$
\sum_{n_1 \ge n_2 \ge \dots \ge n_{l-1} \ge 1}
z^{n_1-1}
\prod_{j=1}^{h-2} \frac{1}{(n_j+p_j)^{u_j}}
\cdot \frac{1}{(n_{h-1}+p_{h-1})^k} \cdot
\prod_{j=h}^{l-1} \frac{1}{(n_j+p_{j+1})^{u_{j+1}}}.
$$
The corresponding parameters in it are
$$
l'=l-1, \quad m'=m+k-u_{h-1}-u_h, \quad
\vec{p}' = (p_1, \dots, p_{h-2}, p_{h-1}, p_{h+1}, \dots, p_l).
$$
Let $P_{\vec{t}}(z)$ be the polynomials in the decomposition
into a linear form in generalized polylogarithms;
then the common denominator of the coefficients of
$P_{\vec{t}}(z)$ divides $D_P^{m'-w(\vec{t})}$.
Since $D_P^{u_{h-1}+u_h-k} A_k \in \Z$,
we have $D_P^{m-w(\vec{t})} (A_k \cdot P_{\vec{t}}(z)) \in \Z[z]$ as required.
The heights of $P_{\vec{t}}(z)$ do not exceed
$$
(l-1)! \cdot (m 2^m)^{l-2} \cdot P^{l-1}.
$$
Consequently, the heights of the polynomials
in the decomposition of sum (\ref{ShiftedCase2Sum1}) do not exceed
\begin{align*}
& \left( \sum_{k=1}^{u_{h-1}} |A_k| + \sum_{k=1}^{u_{h}} |B_k| \right) 
\cdot (l-1)! \cdot (m 2^m)^{l-2} \cdot P^{l-1} \\
\le &
\left( \sum_{k=1}^{u_{h-1}} \binom{u_{h-1}+u_h-k-1}{u_{h-1}-k}
+ \sum_{k=1}^{u_{h}} \binom{u_{h-1}+u_h-k-1}{u_{h}-k} \right) \\
& \quad \times (l-1)! \cdot (m 2^m)^{l-2} \cdot P^{l-1} \\
\le &
( u_{h-1} + u_h ) 2^{u_{h-1}+u_{h}- 2}
\cdot (l-1)! \cdot (m 2^m)^{l-2} \cdot P^{l-1} \\
\le &
m 2^{m-2}
\cdot  (l-1)! \cdot (m 2^m)^{l-2} \cdot P^{l-1} \\
\le &
\frac{1}{2} \cdot (l-1)! \cdot (m 2^m P)^{l-1}.
\end{align*}

Sum (\ref{ShiftedCase2Sum2}) is considered similarly. 
Further, we can apply the induction hypothesis
to the integral $I(p_1, p_2, \dots,$ $p_h-1, \dots, p_l)$.
For all three summands
(\ref{ShiftedCase2Int}), (\ref{ShiftedCase2Sum1}), (\ref{ShiftedCase2Sum2}),
denominators of the coefficients of the polynomial coefficients of
$\Le_{\vec{t}}(z)$ 
in the linear form (\ref{linform}) divide $D_P^{m-w(\vec{t})}$.
The heights of the polynomials $P_{\vec{s}}(z)$
for the initial sum, in case of $\sum_{j=1}^l p_j > 1$, do not exceed
\begin{multline*}
\left( \sum_{j=1}^l p_j - 1 \right) \cdot (l-1)! \cdot (m 2^m P)^{l-1} +
2 \cdot \frac{1}{2} \cdot (l-1)! \cdot (m 2^m P)^{l-1} \\
=
\sum_{j=1}^l p_j \cdot (l-1)! \cdot (m 2^m P)^{l-1}.
\end{multline*}
If $\sum_{j=1}^l p_j = 1$, vectors of the generalized
polylogarithms from the decomposition of
(\ref{ShiftedCase2Sum1}) and (\ref{ShiftedCase2Sum2})
have length less than $l$,
and in the decomposition $I(\inbr{0}_l)$ there is exactly
one polylogarithm of length $l$,
i.e. the sets of the polylogaritms are not intersected
and the estimate on the heights in this case is also valid.

It remains to prove the statement of lemma for the integral
$$
I(p_1, 0, \dots, 0) = \int_{\cube{m}} \frac{
(x_1 x_2 \dots x_{r_1})^{p_1}}
{\prod_{j=1}^l (1-z x_1 x_2 \dots x_{r_j}) }
dx_1 dx_2 \dots dx_{m}.
$$
From the equality
$$
(x_1 x_2 \dots x_{r_1})^{p_1} =
z^{-1} (x_1 x_2 \dots x_{r_1})^{p_1-1} -
z^{-1} (x_1 x_2 \dots x_{r_1})^{p_1-1}(1 - z x_1 x_2 \dots x_{r_1})
$$
it follows that
\begin{align*}
I(p_1, 0, \dots, 0) & = z^{-1} I(p_1-1, 0, \dots, 0)  \\
& \quad
- z^{-1} \int_{\cube{m}} \frac{
(x_1 x_2 \dots x_{r_1})^{p_1-1}}
{\prod_{j=2}^l (1-z x_1 x_2 \dots x_{r_j}) }
dx_1 dx_2 \dots dx_{m} \\
& =z^{-1} I(p_1-1, 0, \dots, 0) \\
& \quad
- z^{-1}
\sum_{n_1 \ge \dots \ge n_{l-1} \ge 1}
z^{n_1-1} \frac{1}{(n_1+p_1-1)^{u_1} n_1^{u_2}}
\prod_{j=2}^{l-1} \frac{1}{n_j^{u_{j+1}}},
\end{align*}
Thus, one can proceed as before in the case $p_h>0$ for $h>1$.
Now the lemma is completely proved.

%% file: Arithm.tex
\def\A{{\mathcal{A}}}

Let us study denominators of the coefficients of the linear form.
We shall use the notion of integer-valued polynomial.
For a polynomial of degree $N$ to be integer-valued
it is sufficient that it possesses integer values
at $N+1$ neighbour integer points (see \cite[Theorem 12.1]{Prasolov}).

Let $\Delta$ be a fixed nonnegative integer.
We say that a rational function $R(x)$ is $\Delta$-normal
if it can be represented as
$$
R(x)=\sum_{\alpha \in \A} \sum_{m=1}^M \frac{A_{m,\alpha}}{(x+\alpha)^m}+
P(x),
$$
where $\A$ is a set of nonnegative integers from a certain
segment $[\alpha_1, \alpha_2]$,
$D_{\Delta}^{M-m} A_{m,\alpha} \in \Z$ and
$D_{\Delta}^{M} P(x)$ is an integer-valued polynomial.

\begin{lemma}
\label{DeltaNorm}
Multiplying 
$\Delta$-normal function by an integer-valued polynomial
of degree $\le \Delta$ remains it $\Delta$-normal.
\end{lemma}
\proof
An integer-valued polynomial $D_{\Delta}^{M} P(x)$,
multiplied by another integer-valued polynomial,
remains the integer-valued.
The statement of the lemma would be proved if we demonstrate it for
$$
R(x)=\frac{A_{m,\alpha}}{(x+\alpha)^m}, \quad
D_{\Delta}^{M-m} A_{m,\alpha} \in \Z.
$$
It is carried out by induction on $m$.
We check firstly the induction base $m=1$.

Let $T(x)$ be an integer-valued polynomial of degree $\le \Delta$
and $\alpha$ be an integer. Then
$$
\frac{T(x)}{x+\alpha}=
\frac{T(-\alpha)}{x+\alpha}+Q(x),
$$
where $Q(x)$ is a polynomial of degree $\le \Delta - 1$
(if $\Delta=0$ it is absent). By the hypothesis, $T(-\alpha)$ is an integer.
Consider $Q(x)$ at the points $x=-\alpha+k$, where $k=1, 2, \dots,
\Delta$:
$$
Q(-\alpha+k)=\frac{T(-\alpha+k)-T(-\alpha)}{k}.
$$ 
Multiplying all these numbers by $D_{\Delta}$ gives
integers, hence $D_{\Delta} Q(x)$ is an integer-valued polynomial.

Thus, if $m=1$,
$$
R(x) T(x) = \frac{A_{1,\alpha} T(-\alpha) }{x+\alpha}+
A_{1,\alpha} Q(x).
$$
In addition,
$$
D_{\Delta}^{M-1} ( A_{1,\alpha} \cdot T(-\alpha) ) =
( D_{\Delta}^{M-1} \cdot A_{1,\alpha} ) T(-\alpha)
\in \Z
$$
and
$$
D_{\Delta}^{M} ( A_{1,\alpha} \cdot Q(x) ) =
( D_{\Delta}^{M-1} \cdot A_{1,\alpha} )\cdot  (D_{\Delta} \cdot Q(x))
$$
is an integer-valued polynomial.

Suppose that $m>1$. Then
$$
R(x) T(x) = \frac{A_{m,\alpha}}{(x+\alpha)^{m-1}} \cdot \frac{T(x)}{x+\alpha}
=\frac{A_{m,\alpha} T(-\alpha)}{(x+\alpha)^m}  + 
\frac{A_{m,\alpha}}{(x+\alpha)^{m-1}} Q(x).
$$
Since
$$
D_{\Delta}^{M-m} \cdot (A_{m,\alpha} T(-\alpha)) =
( D_{\Delta}^{M-m} \cdot A_{m,\alpha} ) T(-\alpha) \in \Z,
$$
the first summand is $\Delta$-normal.
Write the second summand as
$$
\frac{A_{m,\alpha} / D_{\Delta}}{(x+\alpha)^{m-1}} \cdot ( D_{\Delta} Q(x) ).
$$
Since $D_{\Delta}^{M-(m-1)} \cdot (A_{m,\alpha} / D_{\Delta}) \in \Z$
and $D_{\Delta} Q(x)$ is an integer-valued polynomial,
we can apply the induction hypothesis to the latter expression.
The lemma is proved.

Define the {\it index} of
a rational function $R(x) = \frac{P(x)}{Q(x)}$ as
$I(R) = \deg P - \deg Q$.

\begin{lemma}
\label{DeltaSumsUpgrade}
Suppose that for the sum
$$
\F=\sum_{n_1=1}^{\infty} z^{n_1-1} R_1(n_1)
\sum_{n_2=1}^{n_1} R_2(n_2)
\dots
\sum_{n_l=1}^{n_{l-1}} R_l(n_l),
$$
the following inequalities are valid:
\begin{equation}
\label{IndicesSum}
\sum_{j=1}^{j_1} (I(R_j)+1) \le 0, \quad
\sum_{j=j_1}^{j_2} (I(R_j)+1) \le \Delta
\end{equation}
for any $1 \le j_1 \le j_2 \le l$,
and that the functions $R_j$ are $\Delta$-normal.
Then $\F$ is expressed as a finite sum
$\sum_{i} \lambda_i \F_i$, $\lambda_i \in \Q$,
where
$$
\F_i=\sum_{n_1=1}^{\infty} z^{n_1-1} R_{i,1}(n_1)
\sum_{n_2=1}^{n_1} R_{i,2}(n_2)
\dots
\sum_{n_{l(i)}=1}^{n_{l(i)-1}} R_{i,l(i)}(n_{l(i)}),
$$
and $I(R_{i,j}) < 0$ for any $i$, $j$.
In addition, the functions $R_{i,j}$ are $\Delta$-normal
and $D_{\Delta}^{w_i} \lambda_i \in \Z$, where
$$
w_i=\sum_{j=1}^l M_j - \sum_{j=1}^{l(i)} M_{i,j},
$$
$M_{j}$, $M_{i,j}$ is the maximal order of poles of the functions
$R_j$ and $R_{i,j}$.
\end{lemma}
\proof
We proceed by induction on the vector $(l,k)$,
where $k$ is the number of the functions
$R_j$ with $I(R_j) \ge 0$ $(0 \le k < l)$.
Order vectors $(l,k)$ in lexicographic ally.
The induction base $l=1$ is clear, since $I(R_1) \le -1$ in this case
by the hypothesis.
We prove the statement for a vector $(l,k)$
assuming that for smaller vectors it is already proved.
If $k=0$, it is nothing to prove, since $I(R_j) < 0$ for any $j$.
Let $k>0$, i.e. there exists $j$ such that $I(R_j) \ge 0$.
The condition $I(R_1) \le -1$ implies $j>1$.
Expressing $R_j$ as the sum of a polynomial and a proper fraction,
write $\F$ as the sum of two summands.
The summand with with the proper fraction ($\Delta$-normal) has
the number $k$ smaller
by one, hence we can apply the induction base to it.
Now consider the second summand,
in which $R_j(x) = P(x)$ is a polynomial.
From the normality of $R_j$, it follows that the polynomial
$D_{\Delta}^{M_j} P$ is integer-valued
and, in addition, the sum of the maximal orders of poles of
the functions $R_j$ is just smaller by $M_j$ compared with $\F$.

a) If $j=l$ summation of the latter gives
$$
\sum_{n_l=1}^{n_{l-1}} P(n_l) = Q(n_{l-1}),
$$
where $D_{\Delta}^{M_j} Q$ is an integer-valued polynomial of degree
$\deg P + 1$. Thereby, $R_{l-1}$ is multiplied by $Q$.
Thus, compared with the initial sum, the number of summations
is decreased by one.
We can apply the induction hypothesis
to the above sum, multiplied by $D_{\Delta}^{M_j}$,
since the index vector of involved rational functions equals
$(I(R_1)$, \dots, $I(R_{l-2})$, $I(R_{l-1}) + I(R_{l}) + 1)$),
and multiplying by $D_{\Delta}^{M_j} Q(x)$ of
the function $R_{l-1}$ remains it $P$-normal 
by Lemma \ref{DeltaNorm},
because of condition (\ref{IndicesSum}) for $j_1=j_2=j$ we have
$$
\deg Q(x) = \deg P + 1 = I(R_{j}) + 1 \le \Delta.
$$

b) Suppose that $R_j(x)=P(x)$ for $1<j<l$.
Write the initial sum as
$$
\sum_{n_1=1}^{\infty} z^{n_1-1} R_1(n_1)
\sum_{n_2=1}^{n_1} R_2(n_2)
\cdots 
\sum_{n_{j-1}=1}^{n_{j-2}} R_{j-1}(n_{j-1})
\sum_{n_j=1}^{n_{j-1}} P(n_j)
\sum_{n_{j+1}=1}^{n_j} f(n_{j+1}),
$$
where
$$
f(n_{j+1}) = R(n_{j+1})
\sum_{n_{j+2}=1}^{n_{j+1}} R_{j+2}(n_{j+2})
\cdots
\sum_{n_l=1}^{n_{l-1}} R_l(n_l)
$$

We have equalities:
\begin{align*}
& \sum_{n_j=1}^{n_{j-1}} P(n_j)
\sum_{n_{j+1}=1}^{n_j} f(n_{j+1}) \\
& = \sum_{n_j=1}^{n_{j-1}} P(n_j)
\sum_{n_{j+1}=1}^{n_{j-1}} f(n_{j+1})-
\sum_{n_j=1}^{n_{j-1}} P(n_j)
\sum_{n_{j+1}=n_j+1}^{n_{j-1}} f(n_{j+1}) \\
& = Q_1(n_{j-1})
\sum_{n_{j+1}=1}^{n_{j-1}} f(n_{j+1})-
\sum_{n_{j+1}=2}^{n_{j-1}} f(n_{j+1})
\sum_{n_j=1}^{n_{j+1}-1} P(n_j) \\
& = Q_1(n_{j-1})
\sum_{n_{j+1}=1}^{n_{j-1}} f(n_{j+1})-
\sum_{n_{j+1}=1}^{n_{j-1}} Q_2(n_{j+1}) f(n_{j+1});
\end{align*}
in addition, $\deg Q_1 = \deg Q_2 = \deg P + 1$, $Q_2(1)=0$.
Thereby, we express the initial sum as the difference
of sums with a smaller repetition factor,
corresponding to the vectors 
$$( I(R_1), \dots, I(R_{j-1}) + I(R_{j}) + 1, I(R_{j+1}) ,\dots, I(R_{l}) ),$$
$$( I(R_1), \dots, I(R_{j-1}), I(R_{j+1}) + I(R_{j}) + 1,\dots, I(R_{l}) ).$$

The inequality (\ref{IndicesSum}) is also valid for each sum.
Since
$D_{\Delta}^{M_j} P(x)$ is an integer-valued polynomial
and $Q_1(x)$, $Q_2(x)$ are polynomials such that
$Q_1(n)=\sum_{k=1}^n P(k)$, $Q_2(n)=\sum_{k=1}^{n-1} P(k)$
for every integer $n \ge 1$, it follows that $D_{\Delta}^{M_j} Q_1$
and $D_{\Delta}^{M_j} Q_2$ are integer-valued.
Multiplying functions $R_{j-1}$ and $R_{j+1}$
by $D_{\Delta}^{M_j} Q_1(x)$ and $D_{\Delta}^{M_j} Q_2(x)$
remains them $P$-normal by Lemma \ref{DeltaNorm}, because of
$$
\deg Q_1(x) = \deg Q_2(x) = \deg P + 1 = I(R_{j}) + 1 \le \Delta.
$$
The last inequality is due to condition (\ref{IndicesSum}) for $j_1=j_2=j$.
Thus, we can apply the induction hypothesis for each of two sums,
multiplied by $D_{\Delta}^{M_j}$. This completes the proof.

\begin{lemma}
\label{WeakArithmLemma}
Let the parameters $a_i$, $b_i$, $c_j$ be integer,
$b_i > a_i \ge 1$ for $i=1, \dots, m$,
$P=\max_{1 \le i \le m} b_i - 2$, $ q_j = \sum_{i=r_{j-1}+1}^{r_j} (b_i-a_i)$
and let inequalities
$1 \le c_j \le P+1$,
$c_1 + \dots + c_j \le q_1 + \dots + q_j$, $j =1, \dots, l$;
$c_{j_1-1} + \sum_{j=j_1}^{j_2} (c_j-q_j) \le P + 1$,
$1 < j_1 \le j_2 \le l$ be valid.
Suppose that $P_{\vec{s}}$ are the polynomials from the linear form
$S(z)=\sum_{\vec{s}} P_{\vec{s}}(z^{-1}) \Le_{\vec{s}}(z)$.
Then the polynomial
$D_{P}^{m-w(\vec{s})} P_{\vec{s}}(z)$ has integer coefficients.
\end{lemma}
\proof
Using \cite[Lemma 2]{zl5},
represent the integral $S(z)$ as
\begin{multline*}
S(z) =
\frac{\prod_{i=1}^m \Gamma(b_i-a_i)}{\prod_{j=1}^l
\Gamma(c_j)}
\sum_{n_1=1}^{\infty} \sum_{n_2=1}^{n_1} \cdots \sum_{n_l=1}^{n_{l-1}}
 z^{n_1-1}
\\
\times \frac{\prod_{j=1}^l \left[ (n_j-n_{j+1}+1)
(n_j-n_{j+1}+2) \dots (n_j-n_{j+1}+c_j-1) \right]}{\prod_{j=1}^l
\prod_{i=r_{j-1}+1}^{r_j} \left[ (n_j+a_i-1) (n_j+a_i) \cdots
(n_j+b_i-2) \right]},
\end{multline*}
letting $n_{l+1} \equiv 1$.
From the known formula
\begin{multline*}
(x-y+1) \cdots (x-y+n) =
\sum_{k=0}^{n} (-1)^{k} \binom{n}{k} 
(x+k+1) \cdots (x+n) \\ \times
(y + 1) \cdots (y + k - 1)
\end{multline*}
(see, for instance, \cite[Lemma 5]{vasiliev})
it follows that 
\begin{multline*}
(n_j-n_{j+1}+1) \cdots (n_j-n_{j+1}+c_j-1)
=
\sum_{k_j=0}^{c_j-1} (-1)^{k_j} \binom{c_j-1}{k_j} \\
\times
(n_j+k_j+1) (n_j+k_j+2) \cdots
(n_j+c_j-1) \cdot n_{j+1} (n_{j+1} + 1) \cdots (n_{j+1} + k_j - 1).
\end{multline*}
Using this equality for each $j$, we express
$S(z)$ as a linear combination with integer coefficients of sums
(with fixed $k_j$) of type
\begin{multline*}
\sum_{n_1=1}^{\infty} \sum_{n_2=1}^{n_1} \cdots \sum_{n_l=1}^{n_{l-1}}
z^{n_1-1}
\prod_{j=1}^l
p_j^1 (n_j) p_j^2(n_{j+1}) \\
\times
\prod_{i=r_{j-1}+1}^{r_j}
\frac{\Gamma(b_i-a_i)} {(n_j+a_i-1) (n_j+a_i) \cdots (n_j+b_i-2)},
\end{multline*}
where
\begin{align*}
p_j^1(x)&=\frac{(x+k_j+1) (x+k_j+2) \cdots (x+c_j-1)}{(c_j-k_j-1)!}, \\
p_j^2(x)&=\frac{x (x + 1) \cdots (x + k_j - 1)}{k_j!}
\end{align*}
are integer-valued polynomials.
Write the last expression as
$$
\sum_{n_1=1}^{\infty} z^{n_1-1} R_1(n_1)
\sum_{n_2=1}^{n_1} R_2(n_2)
\dots                
\sum_{n_l=1}^{n_{l-1}} R_l(n_l).
$$
where
$$
R_j(x)=
p_j^1(x) p_{j-1}^2(x)
\prod_{i=r_{j-1}+1}^{r_j}
\frac{\Gamma(b_i-a_i)} {(x+a_i-1) (x+a_i) \cdots (x+b_i-2)}.
$$
For $j=1$, define $p_0^2(x) \equiv 1$, $c_0=1$, $k_0=0$.

Since 
$|(b_{i_1}-2)-(a_{i_2}-1)| \le P-(\min\limits_{1 \le i \le m} a_i - 1) \le P$,
the product
$$
\prod_{i=r_{j-1}+1}^{r_j}
\frac{\Gamma(b_i-a_i)} {(x+a_i-1) (x+a_i) \cdots (x+b_i-2)}
$$
is $P$-normal.
Consequently, we can apply Lemma \ref{DeltaNorm} to it
multiplied by $p_j^1(x)$ and $p_{j-1}^2(x)$.
The estimate on polynomial degrees are valid:
$\deg p_j^1 \le c_j-1 \le P$
and $\deg p_{j-1}^2 \le c_{j-1}-1 \le P$.
Thus, $R_j$ is $P$-normal function.

Verify condition (\ref{IndicesSum}) for the function $R_j$:
\begin{align*}
\sum_{j=j_1}^{j_2} (I(R_j)+1) &=
\sum_{j=j_1}^{j_2} (k_{j-1}+(c_j - k_{j} - 1)-q_j+1) \\
& \le c_{j_1-1}-1 + \sum_{j=j_1}^{j_2} (c_j-q_j) \le P.
\end{align*}
The last inequality holds by the hypothesis of the lemma.

Using Lemma \ref{DeltaSumsUpgrade},
we may assume that $I(R_j) < 0$ for any $j$,
and $R_j$ is $P$-normal.
Equivalently, we expressed $S(z)$ as a linear combination
with integer coefficients of sums
$$
\sum_{n_1=1}^{\infty} z^{n_1-1} \frac{A_{\vec{u},\vec{\alpha}}}
{(n_1+\alpha_1)^{u_1}}
\sum_{n_2=1}^{n_1} \frac{1}{(n_2+\alpha_2)^{u_2}}
\dots                
\sum_{n_{l'}=1}^{n_{{l'}-1}} \frac{1}{(n_{l'}+\alpha_{l'})^{u_{l'}}},
$$
where $l' \le l$, $p \le \alpha_j \le P$.
Herewith $D_{P}^{m-w(\vec{u})} A_{\vec{u},\vec{\alpha}} \in \Z$.
Furthermore, for a polynomial $P_{\vec{s}}$ in the
decomposition of the elementary sum
$$
\sum_{n_1=1}^{\infty} z^{n_1-1} \frac{1}
{(n_1+\alpha_1)^{u_1}}
\sum_{n_2=1}^{n_1} \frac{1}{(n_2+\alpha_2)^{u_2}}
\dots                
\sum_{n_{l'}=1}^{n_{{l'}-1}} \frac{1}{(n_{l'}+\alpha_{l'})^{u_{l'}}},
$$
into the linear form $\sum_{\vec{s}} P_{\vec{s}}(z^{-1}) \Le_{\vec{s}}(z)$,
we have the inclusion
$D_{P}^{w(\vec{u})-w(\vec{s})} P_{\vec{s}}(z) \in \Z[z]$ by Lemma \ref{shifted}.
This implies the lemma.

{\bf Remark.} Lemma \ref{WeakArithmLemma} remains valid
if some of $c_j$ are equal to zero.

\begin{theorem}
\label{ArithmTh1}
Let the parameters $a_i$, $b_i$, $c_j$ be integers,
$b_i > a_i \ge 1$ for $i=1, \dots, m$ and $c_j \ge 1$,
$c_1 + \dots + c_j \le q_1 + \dots + q_j$,
where $ q_j = \sum_{i=r_{j-1}+1}^{r_j} (b_i-a_i)$,
$j =1, \dots, l$;
let $d_j$ be nonnegative integers,
satisfying $d_j \le c_j$ for $j=1,\dots,l$ and
$\sum_{k=j}^{l} d_k < a_i$
for $j=1,\dots,l$ and $r_{j-1} < i \le r_j$.
Denote 
$$
\Delta = \max_{1 \le j \le l} \max_{r_{j-1} < i \le r_j}
(b_i - \sum_{k=j}^{l} d_k - 2).
$$
Assume that the following inequalities are valid:
$1 \le c_j \le \Delta + 1$,
$c_{j_1-1} + \sum_{j=j_1}^{j_2} (c_j-q_j) \le \Delta + 1$,
$1 < j_1 \le j_2 \le l$,
and that $P_{\vec{s}}$ are the polynomials from the linear form
$S(z)=\sum_{\vec{s}} P_{\vec{s}}(z^{-1}) \Le_{\vec{s}}(z)$.
Then the polynomial
$D_{\Delta}^{m-w(\vec{s})} P_{\vec{s}}(z)$ has integer coefficients.
\end{theorem}
\proof
Expand the integrand of $S(z)$ using the following equalities:
$$
(x_1 x_2 \cdots x_{r_j})^{d_j}=
\left( \frac{1-(1-z x_1 x_2 \cdots x_{r_j} )}{z} \right)^{d_j},
\quad j=1,\dots,l.
$$
It is possible to do,
since $\sum_{k=j}^l d_k < a_i$ for $j=1,\dots,l$ and $r_{j-1} < i \le r_j$.
This results in a linear combination with integer coefficients 
of expressions of the form
$$
\frac{1}{z^{d_1+d_2+\dots+d_l}}
\int_{\cube{m}} \frac{ \prod_{i=1}^m
x_{i}^{a_{i}'-1} (1-x_{i})^{b_{i}-a_{i}-1} }
{\prod_{j=1}^l (1-z x_1 x_2 \dots x_{r_j})^{c_{j}'} } dx_1 dx_2
\dots dx_{m},
$$
with parameters $c_j'$, satisfying $0 \le c_j' \le c_j$, 
$a_i' = a_i- \sum_{k=j}^l d_k \ge 1$ for $j=1,\dots,l$
and $r_{j-1} < i \le r_j$.
Application of Lemma \ref{WeakArithmLemma}
(in the lemma, $\Delta$ appears as $P$)
to each such integral completes the proof.

\begin{corollary}
\label{cor1}
Let the integral $S(z)$ has parameters
$$
a_i = n+1, \quad b_i=2n+2, \quad c_j=n+1.
$$
Then the polynomial
$D_{n}^{m-w(\vec{s})} P_{\vec{s}}(z)$ has integer coefficients.
\end{corollary}
\proof 
Take  $d_j=0$ for $j=1,\dots,l-1$
and $d_l=n$ in Theorem \ref{ArithmTh1}.
Then $\Delta=n$ and all conditions of the theorem are satisfied.

We apply Corollary \ref{cor1} to the integral
$$
I_{2l+1,n} (z) = \int_{\cube{2l+1}} \frac{ \prod_{i=1}^{2l+1}
x_{i}^n (1-x_{i})^n dx_1 dx_2 \dots
dx_{2l+1}} { \prod_{j=1}^l (1- z x_1 \dots
x_{2j} )^{n+1} (1-z x_1 x_2 \dots x_{2l} x_{2l+1})^{n+1} }.
$$
By \cite[Theorem 6]{zl5}, 
$$
I_{2l+1,n} (z) = 
\sum_{k=0}^{l} P_k (z^{-1}) \Le_{\inbr{2}_k,1} (z)+
\sum_{k=0}^{l-1} T_k (z^{-1}) \Le_{1,\inbr{2}_{k},1} (z)-U(z^{-1}),
$$
($\inbr{a}_k$ means $\{ a, \dots, a \}$,  $k$ times repeated)
where $P_k$, $T_k$, $U$ are polynomials with rational coefficients
and $P_0(1)=0$, $T_k(1)=0$. From Corollary \ref{cor1},
we conclude that these polynomials
multiplied by $D_n^{2l+1}$ have integer coefficients.
Letting $z \to 1-$ and using the equalities
$\Le_{\inbr{2}_k, 1}(1)=2 \zeta(2k+1)$ (see \cite{zl4})
and (\ref{VEqualSOdd}), this proves Vasiliev's conjecture (\ref{VasConj}).

%% file: UpperBound.tex
It is important in many arithmetical applications to have
upper estimates for 
absolute values of the linear form coefficients.
In this section, we study the height of a polynomial
in a linear form in generalized polylogarithms,
that originates from the integral $S(z)$
(see (\ref{linform})).

We start with an estimate for factorial coefficients.

\begin{lemma}
\label{FactorialToPower}
For nonnegative integers $a$ and $b$, the following estimate holds:
$$
\frac{1}{a+b+1} \cdot \frac{(a+b)^{a+b}}{a^a b^b} \le
\frac{(a+b)!}{a! b!} \le
\frac{(a+b)^{a+b}}{a^a b^b}
$$
(if $x=0$, we let $x^x=1$).
\end{lemma}
\proof
If $a=0$ or $b=0$, then both inequalities are valid.
In what follows, suppose that $a$ and $b$ are positive integers.

Consider the Beta-integral
$$
\int_0^1 x^a (1-x)^b dx = B(a+1, b+1) = \frac{a! b!}{(a+b+1)!}.
$$
The function $f(x)=x^a (1-x)^b$, on the segment $[0,1]$,
achives its maximum at the point $x=a/(a+b)$. Hence,
$$
\frac{a! b!}{(a+b+1)!} \le
f \left( \frac{a}{a+b} \right) = \frac{a^a b^b}{(a+b)^{a+b}},
$$
proving the first inequality. Now we prove the second inequality by
induction on the value of $a+b$. The induction base $a=b=1$ is valid.
Introduce the notation
$$
g(a,b)=\frac{(a+b)!}{a! b!}.
$$
Assuming $b>1$, the induction hypothesis yields
$$
g(a,b-1) \le \frac{(a+b-1)^{a+b-1}}{a^a (b-1)^{b-1}}.
$$
From the definition of the function $g$,
$$
\frac{g(a,b)}{g(a,b-1)} = \frac{a+b}{b}.
$$
The function $(1+1/m)^m$ monitonically increases with $m$, hence
$$
\left( 1+\frac{1}{a+b-1} \right)^{a+b-1} \ge
\left( 1+\frac{1}{b-1} \right)^{b-1}.
$$
Write the last inequality as
$$
\frac{(a+b-1)^{a+b-1}}{(b-1)^{b-1}} \le
\frac{(a+b)^{a+b-1}}{b^{b-1}}
$$
Thus,
\begin{align*}
g(a,b)= \frac{a+b}{b} \cdot g(a,b-1) & \le
\frac{a+b}{b} \cdot \frac{(a+b-1)^{a+b-1}}{a^a (b-1)^{b-1}} \\
& \le
\frac{a+b}{b} \cdot \frac{(a+b)^{a+b-1}}{a^a b^{b-1}}=
\frac{(a+b)^{a+b}}{a^a b^b},
\end{align*}
which is the required assertion.

{\bf Remark.} Expression
$$
\frac{(a+b)^{a+b}}{a^a b^b}
$$ 
can be written as
$$
\left( 
\frac{(\alpha+\beta)^{\alpha+\beta}}{\alpha^\alpha \beta^\beta}
\right)^n,
$$
where $\alpha=a/n$, $\beta=b/n$.

By \cite[Lemma 2]{zl5}, the integral $S(z)$
is expressed as
$$
S(z)= \sum_{n_1 \ge n_2 \ge \dots \ge n_l \ge 1 }
R(n_1, n_2, \dots, n_l) z^{n_1-1},
$$
where
\begin{multline*}
R(\zeta_1, \zeta_2, \dots, \zeta_l) =
\frac{\prod_{i=1}^m \Gamma(b_i-a_i)}{\prod_{j=1}^l \Gamma(c_j)} \\*
\times
\frac{\prod_{j=1}^l \left[ (\zeta_j-\zeta_{j+1}+1)
(\zeta_j-\zeta_{j+1}+2) \dots (\zeta_j-\zeta_{j+1}+c_j-1)
\right]}{\prod_{j=1}^l \prod_{i=r_{j-1}+1}^{r_j} \left[
(\zeta_j+a_i-1) (\zeta_j+a_i) \dots (\zeta_j+b_i-2) \right]},
\end{multline*}

To the end of the section, suppose that the parameters
$a_i$, $b_i$, $c_j$ depend linearly on an increasing parameter
$n$, i.e.
$$
a_i=\alpha_i n + \alpha_i', \quad
b_i=\beta_i n + \beta_i', \quad
c_j=\gamma_j n + \gamma_j', \quad
\alpha_i, \beta_i, \gamma_j \in \N, \quad
\alpha_i', \beta_i', \gamma_j' \in \Z.
$$
As before, $q_j = \sum_{i=r_{j-1}+1}^{r_j} (b_i-a_i)$.
We also use notation
$$
p_j = \min\limits_{r_{j-1}+1 \le i \le r_{j}} a_i - 1, \quad
P_j = \max\limits_{r_{j-1}+1 \le i \le r_{j}} b_i - 2,
$$
$$
h_j = \min\limits_{r_{j-1}+1 \le i \le r_{j}} \alpha_i, \quad
H_j = \max\limits_{r_{j-1}+1 \le i \le r_{j}} \beta_i,
$$
$$
\varphi(x, y)=|x+y|^{x+y} \cdot |x|^{-x}.
$$
Here and in what follows,
$|x|^x=1$ if $x=0$, that agrees with the limit value of
$|x|^x$ as $x \to 0$.

\begin{lemma}
\label{UpperBoundLemma}
Let
$c_1 \le q_1$ and $c_{j-1} + c_j \le q_j$ for $j=2, \dots, l$.
Then
$$
R(\zeta_1, \zeta_2, \dots, \zeta_l) = 
\sum_{\vec{s}, \vec{k}} A_{\vec{s}, \vec{k}}
\prod_{j=1}^{l} \frac{1}{(\zeta_j+k_j)^{s_j}},
$$
and 
$$
|A_{\vec{s}, \vec{k}}| \le (F(x_1, \dots, x_l))^{n+ o(n)},
\quad  n \to \infty,
$$
where
$$
x_j=\frac{k_j-p_j}{P_j-p_j} \in [0,1]
$$
and
\begin{align}
\label{UpperBoundFFunction}
F(x_1, \dots, x_l) & =
\prod_{j=1}^{l}
\prod_{i=r_{j-1}+1}^{r_j}
\frac{(\beta_i-\alpha_i)^{\beta_i-\alpha_i}}
{\varphi(\alpha_i-h_j-(H_j-h_j)x_j, \beta_i-\alpha_i)} \nonumber \\
& \quad \times
\prod_{j=1}^{l-1}
\frac{\varphi(h_{j+1}+(H_{j+1}-h_{j+1})x_{j+1}-h_j-(H_j-h_j)x_j,
\gamma_j)}{\gamma_j^{\gamma_j}} \nonumber \\
& \quad \times
\frac{\varphi(h_l+(H_l-h_l)x_l-\gamma_l, \gamma_l)}{\gamma_l^{\gamma_l}}.
\end{align}
\end{lemma}
\proof
Expand the numerator of the function $R$ into the sum of monomials.
Consider any monomial and the corresponding
function $\widehat R(\zeta_1, \zeta_2, \dots, \zeta_l)$.
The degree of the numerator is less than the degree of the denominator
in each variable in $\widehat R$.
Hence, the function $R(\zeta_1, \zeta_2, \dots, \zeta_l)$ 
can be represented as
$$
R(\zeta_1, \zeta_2, \dots, \zeta_l) = 
\sum_{\vec{s}, \vec{k}} A_{\vec{s}, \vec{k}}
\prod_{j=1}^{l} \frac{1}{(\zeta_j+k_j)^{s_j}},
$$

Let $m_j$ be the maximal order of the pole in variable $\zeta_j$.
Cauchy's integral formula for a polycylindrical domain
(see \cite[(1.28)]{Fuks}), applied to partial derivatives of the function
$$
(\zeta_1+k_1)^{m_1} \cdots (\zeta_{l}+k_{l})^{m_{l}}
R(\zeta_1, \zeta_2, \dots, \zeta_{l}),
$$
implies
\begin{multline*}
A_{\vec{s}, \vec{k}} =
\frac{1}{(2 \pi i)^l}
\int_{|\zeta_1+k_1|=\frac{1}{2}} \cdots
\int_{|\zeta_{l}+k_{l}|=\frac{1}{2}}
R(\zeta_1, \zeta_2, \dots, \zeta_{l}) \\
\times
(\zeta_1+k_1)^{s_1-1} \cdots (\zeta_{l}+k_{l})^{s_{l}-1}
d\zeta_1 \cdots d\zeta_{l}.
\end{multline*}
Introduce the function
$$
\Phi(u,v)=|u+v|!^{\sign(u+v)} \cdot |u|!^{-\sign(u)},
$$
defined for integers $u$ and $v$.
The following inequalities holds on the circle $|\zeta_j+k_j|=1/2$:
$$
|(\zeta_j+a_i-1) \cdots (\zeta_j+b_i-2)| \ge
\Phi((a_i-1)-k_j, b_i-a_i-1) e^{o(n)},
$$
$$
|\zeta_{l} (\zeta_{l}+1) \cdots (\zeta_{l}+c_l-2)| \le
\Phi(k_{l}-c_l+1, c_l-1) e^{o(n)},
$$
$$
|(\zeta_j-\zeta_{j+1}+1) \cdots (\zeta_j-\zeta_{j+1}+c_j-1)| \le
\Phi(k_{j+1} - k_j, c_j-1) e^{o(n)}.
$$
We prove only the first inequality (the latter ones are proved similarly).
First, consider the case
$k_j$ lying in the interval $(a_i-1, b_i-2)$.
For $N < k_j$, we have
$$
|\zeta_j+N| = |(k_j - N) - (\zeta_j+k_j)| \ge (k_j-N)-\frac{1}{2},
$$
and, for $N > k_j$,
$$
|\zeta_j+N| = |(N - k_j) + (\zeta_j+k_j)| \ge (N - k_j)-\frac{1}{2}.
$$
Consequently (Set a product to be 1 if upper limit greater than lower),
\begin{align*}
|(\zeta_j+a_i-1) \cdots (\zeta_j+b_i-2)| & \ge
\prod_{N=a_i-1}^{k_j-1} (k_j-N-\frac{1}{2}) \cdot \frac{1}{2} 
\cdot \prod_{N=k_j+1}^{b_i-2} (N - k_j - \frac{1}{2}) \\
& \ge \frac{1}{8}
\prod_{N=a_i-1}^{k_j-2} (k_j-N-1) \cdot 
\prod_{N=k_j+2}^{b_i-2} (N - k_j - 1) \\
& = \frac{(k_j - (a_i -1))! ((b_i-2)-k_j)!}{8 (k_j-(a_i-1)) ((b_i-2)-k_j)} \\
& = \Phi((a_i-1)-k_j, b_i-a_i-1) e^{o(n)}.
\end{align*}
Now consider the case $k_j<a_i-1$:
\begin{align*}
|(\zeta_j+a_i-1) \cdots (\zeta_j+b_i-2)| & \ge
\prod_{N=a_i-1}^{b_i-2} (N - k_j - \frac{1}{2}) 
\ge \frac{1}{2} \prod_{N=a_i}^{b_i-2} (N - k_j - 1) \\
& =\frac{a_i-1-k_j}{2(b_i-2-k_j)} \cdot
\frac{((b_i-2)-k_j)!}{((a_i-1)-k_j)!} \\
& = \Phi((a_i-1)-k_j, b_i-a_i-1) e^{o(n)}.
\end{align*}
The case $k_j > b_i - 2$ is considered similarly.
Cases $k_j=a_i-1$ or $k_j=b_i-2$ can be verified by direct
substitution.

Thus,
\begin{multline}
\label{CoefBoundN}
|A_{\vec{s}, \vec{k}}| \le 
\prod_{j=1}^{l}
\prod_{i=r_{j-1}+1}^{r_j} \frac{(b_i-a_i-1)!}{\Phi((a_i-1)-k_j, b_i-a_i-1)} \\
\times
\prod_{j=1}^{l-1} \frac{\Phi(k_{j+1}-k_j, c_j - 1)}{(c_j - 1)!} \cdot
\frac{\Phi(k_l-c_l+1, c_l - 1)}{(c_l - 1)!} \cdot
e^{o(n)}.
\end{multline}
Make the change $k_j=p_j+ (P_j-p_j) x_j$, $j=1,\dots,l$, $x_j \in [0,1]$.
From Lemma \ref{FactorialToPower} and estimate (\ref{CoefBoundN})
it follows that
$$
|A_{\vec{s}, \vec{k}}| \le  
(F(x_1, \dots, x_l))^n e^{o(n)},
$$
where $F(x_1, \dots, x_l)$ is the function, specified in the lemma statement.
                                                            
\begin{theorem}
\label{UpperBoundTh}
Let
$c_1 \le q_1$ and $c_{j-1} + c_j \le q_j$ for $j=2, \dots, l$.
Then the heights of all polynomials in the linear form
$S(z)=\sum_{\vec{s}} P_{\vec{s}}(z^{-1}) \Le_{\vec{s}}(z)$ 
do not exceed $M^{n+o(n)}$, as $n \to \infty$,
where $M$ is the maximum of function (\ref{UpperBoundFFunction})
on the cube $\cube{l}$.
\end{theorem}
\proof
By Lemma \ref{UpperBoundLemma}, we have the equality:
\begin{align*}
S(z)
& =\sum_{n_1 \ge n_2 \ge \dots \ge n_l \ge 1 }
R(n_1, n_2, \dots, n_l) z^{n_1-1} \\
& =
\sum_{\vec{s}, \vec{k}} A_{\vec{s}, \vec{k}}
\sum_{n_1 \ge n_2 \ge \dots \ge n_l \ge 1 } z^{n_1-1}
\prod_{j=1}^{l} \frac{1}{(\zeta_j+k_j)^{s_j}}.
\end{align*}
Since $s_j \le m$ and $b_i - a_i \le C n$, then the number of
summands in the external sum does not exceed
$(m \cdot Cn)^l = e^{o(n)}$.
Furthermore, consider the decomposition of the elementary sum
$$
\sum_{n_1 \ge n_2 \ge \dots \ge n_{l} \ge 1 } z^{n_1-1}
\prod_{j=1}^{l} \frac{1}{(n_j+k_j)^{s_j}}
$$
into a linear form in polylogarithms
for a fixed $\vec{s}$ and $\vec{k}$.
The heights of all polynomials in it is $e^{o(n)}$ by Lemma \ref{shifted}.
From Lemma \ref{UpperBoundLemma} it follows that
$|A_{\vec{s}, \vec{k}}| \le M^{n+o(n)}$.
This implies the statement of theorem.

%% file: article.bbl
\begin{thebibliography}{99}

\bibitem{zl5}
{\namefont Zlobin~S.A.,}
{\titlefont Expanding multiple integrals into linear forms} //
Matem. Zametki. 2005. V. 77. \nomer 5. P. 683--706.

\bibitem{vasiliev}
{\namefont Vasilyev~D.V.,}
{\titlefont On small linear forms for the values
of the Riemann zeta-function at odd integers} // Preprint \nomer 1
(558). Minsk: Nat. Acad. Sci. Belarus, Institute Math., 2001.

\bibitem{beukers}
{\namefont Beukers F.,}
{\titlefont A note on the irrationality of $\zeta(2)$ and $\zeta(3)$} //
Bull. London Math Society. 1979. V. 11. \nomer 3. P. 268--272.

\bibitem{zudilin}
{\namefont Zudilin~W.V.,}
{\titlefont Well-poised 
hypergeometric series and multiple integrals} //
Uspehi Matem. Nauk. 2002. V. 57. \nomer 4. P. 177--178.

\bibitem{rivkrat}
{\namefont Krattenthaler C., Rivoal T.,}
{\titlefont Hyperg\'eom\'etrie et fonction z\^eta de Riemann} //
Preprint (December 2004), submitted for publication; //
http://arxiv.org/abs/math/0311114.

\bibitem{zl1}
{\namefont Zlobin~S.A,}
{\titlefont
Integrals expressible as linear forms in generalized polylogarithms} //
Matem. Zametki. 2002. V. 71. \nomer 5. P. 782--787.

\bibitem{zl4}
{\namefont Zlobin~S.A.,}
{\titlefont Generating functions for multiple zeta values} //
Vestnik MGU. Ser. 1. Matem., Mekh. 2005. \nomer 2. P. 55--59.

\bibitem{Prasolov}
{\namefont Prasolov~V.V.,}
{\titlefont Polynomials} //
M.: MCCME, 1999.

\bibitem{Fuks}
{\namefont Fuchs~B.A.,}
{\titlefont
Introduction to the theory of analytic functions of several complex variables}
// M.: Gos. Izd. Fiz.-mat. Lit., 1962.

\end{thebibliography}
